\documentclass[twoside,bezier,reqno]{amsart}

\usepackage{epsfig}
\usepackage{amsmath,amsthm,amsfonts,amssymb,amscd,euscript}

\input{epsf}
\newcommand{\filebegin}{\begin{document}}
\newcommand{\fileend}{\end{document}}
\makeatother
\def\thefootnote{}
\newcommand{\lo}{\longrightarrow}
\newcommand{\NMM}{\hspace*{2mm}}

\renewcommand{\baselinestretch}{1.1}

\input{epsf}
\renewcommand{\baselinestretch}{1.1}
\def\n{\noindent}%

\numberwithin{equation}{section}
\def\mapdown#1{\Big\downarrow\rlap
{$\vcenter{\hbox{$\scriptstyle#1$}}$}}
\newtheorem{theorem}{Theorem}[section]
\newtheorem{lemma}[theorem]{Lemma}
\newtheorem{proposition}[theorem]{Proposition}
\newtheorem{corollary}[theorem]{Corollary}

\theoremstyle{definition}
\newtheorem{definition}[theorem]{Definition}
\newtheorem{example}[theorem]{\sc Example}
\newtheorem{xca}[theorem]{Exercise}
\theoremstyle{remark}
\newtheorem{remark}[theorem]{Remark}
\newcommand{\Sym}{{\rm Sym}}
\newcommand{\Fix}{{\rm Fix}}

\begin{document}

\vspace*{2cm}
\begin{center}
{\bf\large Fixity of elusive groups and the polycirculant conjecture}
 \\[0.5cm]
{Majid Arezoomand  \\[2mm]
University of Larestan, Larestan, 74317-16137, Iran\\[2mm]
{\tt E-mail:arezoomand@lar.ac.ir}}
 \\[2mm]
\end{center}%
\vspace*{0.5cm}
\begin{quotation}
\noindent
{\footnotesize
{\sc Abstract.}
Let $G\leq\Sym(\Omega)$ be transitive. Then $G$ is called \textit{elusive} on $\Omega$  if it has no fixed point 
free element of prime order. The \textit{$2$-closure} of $G$, denoted by $G^{(2),\Omega}$, is the largest subgroup
of $\Sym(\Omega)$ whose orbits on $\Omega\times\Omega$ are the same orbits of $G$. 
$G$ is called $2$-closed on $\Omega$ if $G=G^{(2),\Omega}$. 
The \textit{polycirculant conjecture} states that there is no $2$-closed elusive group. In this paper, we study the \textit{fixity} of elusive groups, where the fixity of $G$ is 
the maximal number of fixed points of a non-trivial element of $G$. In particular, we prove that
there is no $2$-closed elusive solvable group of fixity at most $5$, a partial answer to the polycirculant conjecture.}
\end{quotation}

\n \textbf{2000 Mathematics subject classification: } 20B25, 05E18.

\markboth 
{Majid Arezoomand}
 {Fixity of elusive groups and the polycirculant conjecture}



\section{Introduction}
The concept of fixity was introduced by Rosen in 1980 to study of permutation groups of prime power order. Let $G$ be a permutation group on a set $\Omega$. 
We say that  $G$ has \textit{fixity} $f=f(G)$ if no non-trivial element of $G$ fixes more than $f$ letters, and there is a non-trivial element of $G$ fixing exactly $f$ letter \cite{SS}.
 
Let $G$ be a transitive permutation group. Then $f(G)=0$ if and only if $G$ is
regular. Furthermore, $f(G)=1$ if and only if
 $G$ is a Frobenius group. For some applications of the fixity of permutation groups, we refer the reader to \cite{LS,MW,SS}.

A permutation $g$ of a set $\Omega$ is said to be a {\it derangement} if it has no fixed-point on $\Omega$, equivalently for each $\alpha\in\Omega$, $\alpha^g\neq\alpha$. A transitive permutation group  $G\leq\Sym(\Omega)$ is called \textit{elusive} on $\Omega$  if it has no fixed point 
free element of prime order. Let $G$ be a permutation group on a set $\Omega$. The \textit{$2$-closure} of $G$, denoted by $G^{(2),\Omega}$, is the largest subgroup
of $\Sym(\Omega)$ whose orbits on $\Omega\times\Omega$ are the same orbits of $G$. Then 
$G$ is called $2$-closed on $\Omega$ if $G=G^{(2),\Omega}$. 

In 1981, Maru\v{s}i\v{c} asked
whether there exists a vertex-transitive digraph without a non-identity automorphism having  all of its orbits of the same length
\cite[Problem 2.4]{Marusic}. In 1988, independently, the above problem was again proposed by Jordan~\cite{Jordan}.
In the 15th British  combinatorial conference, in 1995, Klin proposed a more general question in the context of $2$-closed groups \cite{Klin}:
``Is there a $2$-closed
transitive  permutation group containing no fixed-point-free element of prime order?''  After decades, not only has a positive answer to Klin's question not been found, based on the evidence, mathematicians conjecture that there is no positive answer to this question. This conjecture is known as {\it polycirculant conjecture}. Equivalently, the polycirculant conjecture states that
no $2$-closed transitive permutation group is elusive. The conjecture is still open. We refer the reader
to \cite{AAS} to a survey on the recent results and future directions of the polycirculant conjecture.

In this paper, we study the fixity of elusive groups. Then we prove that every transitive $2$-closed solvable permutation group of fixity at most $5$ confirms the polycirculant conjecture.

\section{Main Results}
First we collect some notations we need later. Let $G$ be a finite group and $\Omega$ be a non-empty set. We denote the center of $G$ and the set of all prime divisors the order of $G$ by $Z(G)$ and $\pi(G)$, respectively. Also $\Sym(\Omega)$ denotes the group of all permutations on $\Omega$. Let $G$ acts on $\Omega$ and $\alpha\in\Omega$. We denote the stabilizer of $\alpha$ in $G$ and the orbit of $G$ containing $\alpha$ by $G_\alpha$ and $\alpha^G$, respectively. Also $\Fix_\Omega(G)$ denotes the fixed points of $G$ on $\Omega$, the set of all elements of $\Omega$ which fixes by all elements of $G$. Furthermore, $G^\Omega$ denotes the homomorphic image of the action of $G$ on $\Omega$ which is a subgroup of $\Sym(\Omega)$. Recall that $G$ is called a permutation group on $\Omega$ if $G$ is isomorphic to $G^\Omega$. For the notations and terminology not defined here, we refer the reader to \cite{Dixon}.

Let $G$ be a permutation group on a set $\Omega$. Recall that fixity of $G$ is 
the maximal number of fixed points of a non-trivial element of $G$ on $\Omega$. 
We start with the following lemma:

\begin{lemma}\label{main}
Let $G$ be a finite permutation group on $\Omega$ of fixity $f\geq 2$, $p>f$ a prime and $\alpha\in\Omega$. 
If $p\in\pi(G_\alpha)$ then $G_\alpha$ contains a Sylow $p$-subgroup of $G$. In particular, if $G$ is transitive and contains 
a non-trivial normal $p$-subgroup then $p\notin\pi(G_\alpha)$.
\end{lemma}
\proof
Let $p$ be a prime divisor of order of $G_\alpha$. Then there exists a Sylow $p$-subgroup $P$ of $G$ such that $P_\alpha\neq 1$. We claim that $Z(P)\leq G_\alpha$. Suppose towards a contradiction, that our claim is not true. Then there exists $g\in Z(P)$ such that $g\notin G_\alpha$.  Hence 
$\alpha,\alpha^g,\ldots,\alpha^{g^{p-1}}$ are distinct elements of $\Fix_\Omega(P_\alpha)$,  because
if $\alpha^{g^i}=\alpha^{g^j}$, for
some $0\leq i<j\leq p-1$, then $g^{j-i}\in G_\alpha$ and so $(j-i,p)=1$ implies that $g\in G_\alpha$, a contradiction. So $P_\alpha\neq 1$ implies that
$p\leq f$, a contradiction. This proves that $Z(P)\leq G_\alpha$.

Now we prove that $P\leq G_\alpha$. Suppose, by a contrary, that $P\nleq G_\alpha$. Then $x\notin G_\alpha$ for some $x\in P$. Hence, by a similar argument to the above paragraph, $\alpha,\alpha^x,\ldots,\alpha^{x^{p-1}}$ are distinct. Now $Z(P)\leq G_\alpha$ implies that these elements are all in $\Fix_\Omega(Z(P))$. Let $y$ be a non-trivial element of $Z(P)$. Then
\[p\leq |\Fix_\Omega(Z(P))|=|\bigcap_{x\in Z(P)}\Fix_\Omega(x)|\leq |\Fix_\Omega(y)|\leq f,\]
which is a contradiction.

Finally, suppose that $G$ is transitive and $N$ is a non-trivial normal $p$-subgroup of $G$. Then is contained in every Sylow $p$-subgroup of $G$. Now if $p$ is a divisor of $|G_\alpha|$, then by the above paragraph, $N\leq G_\alpha$. Since $G$ is transitive, this implies that $N$ is contained in any point-stabilizer of $G$ which implies that $N=1$, because $G$ is a permutation group. This completes the proof.
\qed

\begin{corollary}
Let $G$ be an elusive group on $\Omega$ with fixity $f$. If $G$ contains a non-trivial 
normal $p$-subgroup then $p\leq f$.
\end{corollary}
\proof Suppose, towards a contradiction, that $p>f$. Then by Lemma \ref{main}, $p\notin\pi(G_\alpha)$, where $\alpha\in\Omega$. On the other hand, by \cite[Lemma 2.1]{GMPV}, $\pi(G)=\pi(G_\alpha)$ which is a contradiction.
\qed

\begin{corollary}\label{geq3}
Let $G$ be an elusive group on $\Omega$ of fixity $f$. Then $f\geq 3$.
\end{corollary}
\proof
Suppose, towards a contradiction, that $f\leq 2$. Since $f=0$ if and only if $G$ is regular on $\Omega$, and since regular groups are not elusive, we have $f=1,2$.
On the other hand, $f=1$ if and only if $G$ is a Frobenius group with Frobenius complement $G_\alpha$ on $\Omega$, where $\alpha\in\Omega$. Since the Frobenius kernel
of any finite Frobenius group is a regular subgroup, we conclude that $f=2$. 

Let $\alpha\in\Omega$. Then, by \cite[Lemma 2.1]{GMPV}, $\pi(G)=\pi(G_\alpha)$. Hence, by Lemma \ref{main}, $|\Omega|=2^k$,
for some $k\geq 1$, which contradicts \cite[Lemma 2.6]{AAS}.
\qed

\begin{lemma}\label{2,3} 
Let $G$ be an elusive group of fixity $f\geq 3$ on a finite set $\Omega$. Let $|\Omega|=p_1^{n_1}\ldots p_k^{n_k}$,
where $n_i\geq 1$ and $k\geq 2$. Then 
\begin{itemize}
\item[$(i)$] for each $1\leq i\leq k$, $p_i\leq f$,
\item[$(ii)$] if  $p\in\{p_1,\ldots,p_k\}$ then every $p$-element in $G$ is either fixed-point free or fixes $mp$ points, where $1\leq m\leq f/p$.
\end{itemize}
\end{lemma}
\proof
$(i)$ Let $\alpha\in\Omega$ and $p$ be a prime divisor of $|\Omega|=|G:G_\alpha|$. By Lemma \ref{main}, if $p>f$ then $G_\alpha$ contains a Sylow $p$-subgroup of $G$. Hence $p$ can not divide $|\Omega|$, a contradiction.

$(ii)$ Let $p\in\{p_1,\ldots,p_k\}$ and $x\in G$ be a $p$-element. Then $x$ fixes at most $f$ points.
Since the length of any orbit of $\langle x\rangle$ on $\Omega$ is a power of $p$, by $(i)$, $|\Omega|=|\Fix_\Omega(\langle x\rangle)|+lp$,
 for some positive integer $l$. Hence $p$ divides $|\Fix_\Omega(\langle x\rangle)|$. On the other hand, $\Fix_\Omega(\langle x\rangle)=\Fix_\Omega(x)$ which implies that
 $x$ is either fixed-point free or fixes $mp$ points, where $1\leq m\leq f/p$. 
\qed

\begin{corollary} Let $G$ be an elusive group on $\Omega$ with fixity $f$. If $|\Omega|$ is odd, then $f\geq 5$. In particular, if $|G|$ is odd
then $f\geq 5$.
\end{corollary}
\proof
Since $|\Omega|$ is a divisor of $|G|$, the second part is a direct consequence of the first part. Suppose, by contrary, that $f\leq 4$. Then by Corollary \ref{geq3}, $f\in\{3,4\}$. Now by Lemma \ref{2,3}, $|\Omega|=2^m3^n$, where $m,n\geq 0$ are integers. On the other hand, by \cite[Lemma 2.6]{AAS}, $m,n\neq 0$. Let $x$ be an  element of order $3$ of $G$. Then, by Lemma \ref{2,3}, $x$ is fixed-point free, which is a contradiction.
\qed

For a finite group $1\neq G$, we denote the smallest prime divisor of $|G|$ by $p_G$. Then we have the following lemma.

\begin{lemma}\label{min}
Let $G$ be an  elusive group on a finite set $\Omega$ with fixity $f$. Then
\[|\Omega|\leq f\min\{\frac{|H|-1}{p_H-1}\mid 1\neq H\unlhd G\}.\]
\end{lemma}
\proof
Let $H$ be a non-trivial normal subgroup of $G$.  Let $\Omega_1,\ldots,\Omega_m$ be all $H$-orbits on $\Omega$. By \cite[Theorem 1.6. A]{Dixon}, 
for $i\neq j$, $H^{\Omega_i}$ is permutation isomorphic to $H^{\Omega_j}$. If $H$ acts regularly on one orbit, then it acts regularly on all of its orbits and so every element of prime order in $H$ must be fixed-point free which is a contradiction. Hence $H$ does not act regularly on its orbits. Thus
 by \cite[Lemma 2.6]{SS}, $|\Omega|\leq f(|H|-1)/(p_H-1)$, which completes the proof.
\qed

\begin{lemma}\label{abnorm}
Let $G$ be an elusive group on a finite set $\Omega$ with fixity $f$. If $G$ has a non-trivial normal abelian subgroup $N$ 
and $p$ is the smallest prime divisor of $|N|$, then
\begin{itemize}
\item[$(1)$] $|N|\leq pf$,
\item[$(2)$] $|\Omega|\leq \frac{f(pf-1)}{p-1}\leq f(2f-1)$. 
\end{itemize}
\end{lemma}
\proof It is obvious that the fixity is at most $f$. On the other hand, by the proof of Lemma \ref{min}, $N$ has no regular orbit on $\Omega$. Hence, by \cite[Lemma 2.7]{SS}, $|N|\leq pf$. The second part follows from Lemma \ref{min}.
\qed

\begin{corollary} 
Let $G$ be a $2$-closed elusive solvable group of fixity $f$. Then $f\geq 6$.
\end{corollary}
\proof
Let $G$ be a $2$-closed elusive group on $\Omega$. Suppose that $N$ is a minimal normal subgroup of $G$. Then $N\cong\Bbb Z_p^k$ for some $k\geq 1$, where $p$ is a prime. 

Suppose, towards a contradiction, that $f\leq 5$. Then, by Lemma \ref{abnorm}, $|\Omega|\leq 45$ which contradicts \cite[Proposition 6.1]{HR}.
\qed

\begin{corollary}\label{abelian} Let $G$ be an elusive group of fixity $f$ on $\Omega$ and $N\neq 1$ be an abelian
normal subgroup of $G$ of order $p_1^{n_1}p_2^{n_2}\ldots p_k^{n_k}$, where $p_1<p_2<\ldots<p_k$ are primes,
$n_i\geq 1$ and $k\geq 1$. Then
\begin{itemize}
\item[$(1)$] If $k=1$ then $n_1\neq 1$. Also $(n_1,\ldots,n_k)\neq (1,\ldots,1)$.
\item[$(2)$] $p_1^{n_1+\ldots+n_k-1}\leq f$. In particular, $p_1^{k-1}\leq f$.
\item[$(3)$] If $f=3$, then $p_1=2$ or $3$ and $N\cong\Bbb Z_2^2$ or $\Bbb Z_3^2$, respectively.
\item[$(4)$] If $f=4$ then $p_1=2$ or $3$. In the first case $N$ is isomorphic to one of the groups $\Bbb Z_2^2$, $\Bbb Z_2\times \Bbb Z_p^2$ or $\Bbb Z_2^2\times\Bbb Z_p$, where $p$ is an odd prime. In the later case, $N\cong\Bbb Z_3^2$.
\end{itemize}
\end{corollary}
\proof
$(1)$ If $|N|=p_1$ or $p_1p_2\ldots p_k$ then there exists a non-trivial normal cyclic subgroup of $G$ which contradicts \cite[Lemma 2.20]{AAS}.

$(2)$ It is an immediate consequence of Lemma \ref{abnorm}.

$(3)$ Since $G$ is elusive, $N$ is not cyclic by \cite[Lemma 2.20]{AAS}. Now $(3)$ is a consequence of $(1)$ and $(2)$.

$(4)$ Let $f=4$. Then $(2)$ implies that $p_1=2$ or $p_1=3$. If $k=1$ then by Lemma \ref{abnorm}, $N\cong\Bbb Z_2^2$. If $k\geq 2$ the the result follow from $(1)$ and $(2)$.
\qed

\begin{corollary} Let $G$ be a transitive $2$-closed  permutation group on a set $\Omega$ of fixity $4$ and $N\neq 1$ be a normal $p$-subgroup
of $G$. Then $G$ admits fixed-point free element.
\end{corollary}
\proof
Suppose, towards a contradiction, that $G$
is elusive. Then, by Corollary \ref{abelian}, 
we have $N\cong\Bbb Z_2^2$ or $N\cong\Bbb Z_3^2$.  Let $\alpha^N$ be an orbit of $N$ on $\Omega$. Then $|\alpha^N|\in\{1,p,p^2\}$, where $p\in\{2,3\}$. If $|\alpha^N|=1$, then $N=N_\alpha\leq G_\alpha$ which implies that $N=1$, a contradiction. If $|\alpha^N|=p^2$ then $N_\alpha=1$ which contradicts \cite[Lemma 2.7]{AAS}. Hence $|\alpha^N|=p$ which contradicts \cite[Theorem 2.11]{AAS}.
\qed

\section*{Acknowledgments}
The author gratefully appreciate an anonymous referee for constructive comments and recommendations which
definitely helped to improve the readability and quality of the article.

\providecommand{\bysame}{\leavevmode\hbox
to3em{\hrulefill}\thinspace}


\end{document}